\documentclass[review]{elsarticle}
\usepackage{graphicx}
\usepackage{epstopdf}
\usepackage{lineno,hyperref}
 \usepackage{geometry}
\usepackage{epsfig}
\usepackage{lineno,hyperref}
\usepackage{bookmark}
\usepackage{multirow}
 \usepackage{amsfonts}
\usepackage{multirow}
 \usepackage{algorithm}
 \usepackage[figuresright]{rotating}
 \usepackage{amscd}
\usepackage{mathrsfs}
\usepackage{booktabs,longtable}
\usepackage{indentfirst}
\usepackage{subfigure}
\usepackage{amstext}
\usepackage{amssymb}
\usepackage{fancyhdr}
\usepackage{amsmath}
\usepackage{color}
\modulolinenumbers[5]
\biboptions{sort&compress}
\newtheorem{thm}{Theorem}
\newtheorem{lem}[thm]{Lemma}
\newdefinition{rmk}{Remark}
\newproof{pf}{Proof}
\newproof{pot}{Proof of Theorem \ref{thm2}}
\newtheorem{alg}{Algorithm}
\newdefinition{exa}{Example}

\journal{Journal of \LaTeX\ Templates}

\begin{document}

\begin{frontmatter}
\title{Preconvergence of the randomized extended Kaczmarz method
\tnoteref{mytitlenote}}
\tnotetext[mytitlenote]{The work is supported by the National Natural Science Foundation of China (No. 11671060) and the Natural Science Foundation Project of CQ CSTC (No. cstc2019jcyj-msxmX0267)}

\author{Yanjun Zhang, \ Hanyu Li\corref{mycor}}
\cortext[mycor]{Corresponding author. E-mail addresses: yjzhang@cqu.edu.cn;\ lihy.hy@gmail.com or hyli@cqu.edu.cn.}

\address{College of Mathematics and Statistics, Chongqing University, Chongqing 401331, P.R. China}
\begin{abstract}
In this paper, we analyze the convergence behavior of the randomized extended Kaczmarz (REK) method 
for all types of linear systems (consistent or inconsistent, overdetermined or underdetermined, full-rank or rank-deficient). The analysis shows that the larger the singular value of 
$A$ is, the faster the error decays in the corresponding right singular vector space, and as $k\rightarrow\infty$, $x_{k}-x_{\star}$ tends to the right singular vector corresponding to the smallest singular value of $A$, where $x_{k}$ is the $k$th approximation of the REK method and $x_{\star}$ is the minimum $\ell_2 $-norm least squares solution.
These results 
explain the phenomenon found in 
the extensive numerical experiments appearing in the literature that the REK method seems to converge faster in the beginning. A simple numerical example is provided to confirm the above findings.
\end{abstract}


\end{frontmatter}


\section{Introduction}
The Kaczmarz method \cite{kaczmarz1} is a popular iterative method for solving the following linear system
\begin{equation}
\label{1}
Ax=b,
\end{equation}
where $A\in R^{m\times n}$, $b\in R^{m}$, and $x$ is the $n$-dimensional unknown vector, 
and has found a wide range of applications in many fields, such as 
medical scanner \cite{hounsfield1973computerized}, digital signal processing \cite{byrne2003unified, lorenz2014sparse}, distributed computing \cite{elble2010gpu}, computer tomography \cite{censor1988parallel}, image reconstruction \cite{eggermont1981iterative, herman1993algebraic, popa2004kaczmarz}, etc. At each step, the Kaczmarz method orthogonally projects the current estimate onto one hyperplane defined by the $i$th constraint of the system. The convergence of the method is not difficult to show, but the theoretical analysis of its convergence rate is a big challenge.

If the system (\ref{1}) is consistent, Strohmer and Vershynin \cite{Strohmer2009} proved the linear convergence of the randomized Kaczmarz (RK) method for overdetermined full rank linear system. In fact, the RK method has the same convergence property regardless of whether the system is overdetermined or underdetermined, full rank or rank deficient; see \cite{Completion2015, gower2015stochastic} for more details. 
Now, the RK method 
has been extended to solve various problems including linear constraint problem \cite{leventhal2010randomized}, ridge regression problem \cite{hefny2017rows, Liu2019}, linear feasibility problem \cite{de2017sampling}, generalized phase retrieval problem \cite{wei2015solving}, and inverse problem \cite{li2018averaged}, and 
has many 
variants \cite{Eldar2011, needell2014paved, nutini2016convergence, Bai2018, Wu2020, Chen2020}.

If the system (\ref{1}) is inconsistent, 
it holds that $b=Ax_{\star}+z$, where 
$x_{\star}=A^{\dag}b$ is the minimum $\ell_2 $-norm least squares solution with $A^{\dag}$ denoting the Moore-Penrose pseudoinverse of the matrix $A$ and 
$z$ is a nonzero vector belonging to the null space of $A^{T}$.
In this case, Needell \cite{Needell2010} proved that the RK method does not converge to $x_{\star}$. To resolve this convergence problem, Zouzias and Freris \cite{Completion2013} proposed the randomized extended Kaczmarz (REK) method, which essentially uses the RK method twice in each iteration \cite{liu2016accelerated, Dukui2019, du2021randomized}. The REK method can also be considered as a randomized variant of the extended Kaczmarz method proposed by Popa \cite{popa1995least, popa1998extensions}. Later, many variants of the REK method were proposed to accelerate the convergence; see for example \cite{needell2015randomized, xiang2017accelerated, du2020randomized} and references therein.

In 2017, Jiao, Jin and Lu \cite{jiao2017preasymptotic} analyzed the preasymptotic convergence of the RK method. By decomposing a space into two orthogonal subspaces, i.e., the low right singular vectors subspaces (corresponding to the large singular values) and the high right singular vectors subspaces (corresponding to the small singular values), they showed that during initial iterations the error in the low right singular vectors subspaces decays faster than that in the high right singular vectors subspaces. Recently, Steinerberger \cite{steinerberger2021randomized} made a more detailed analysis of the convergence property of the RK method for overdetermined full rank linear system. The author showed that the right singular vectors of the matrix $A$ describe the directions of distinguished dynamics and the RK method converges along small right singular vectors.

In this paper, we are going to take analysis on the convergence property of the REK method for all types of linear systems (consistent or inconsistent, overdetermined or underdetermined, full-rank or rank-deficient). We show that the sequence $\{x_{k}\}_{k=1}^{\infty}$ generated by the REK method converge to the minimum $\ell_2 $-norm least squares solution $x_{\star}$ with different decay rates in different right singular vectors spaces, and as $k\rightarrow\infty$, $x_{k}-x_{\star}$ finally tends to the right singular vector corresponding to the smallest singular value of $A$. 

The rest of this paper is organized as follows. We first introduce some notations and preliminaries in Section \ref{sec2} and then present our main results in Section \ref{sec3}. A simple numerical experiment is given in Section \ref{sec4}. 

\section{Notations and preliminaries }\label{sec2}
Throughout the paper, for a matrix $A$, $A^T$, $A^{(i)}$, $A_{(j)}$,
$\sigma_i(A)$, $\sigma_r(A)$, $\|A\|_F$, and $\mathcal{R}(A)$ denote its transpose, $i$th row (or $i$th entry in the case of a vector), $j$th column, $i$th singular value, smallest nonzero singular value, Frobenius norm, and column space, respectively. For any integer $m\geq1$, let $[m]:=\{1, 2, 3, ..., m\}$. In addition, we denote 
the expectation of any random variable $\xi$ by $\mathbb{E} [\xi]$.

We list the REK method presented 
in \cite{ Dukui2019} in Algorithm \ref{alg1}, which is a slight variant of the original REK method \cite{Completion2013}. From the algorithm we find that, in each iteration, $z_k$ is the $k$th approximation of the RK method applied to $A^Tz=0$ with initial guess $z_0$, and $x_k$ is a one-step RK update for the linear system $Ax=b-z_{k}$ from $x_{k-1}$.
\begin{alg}
\label{alg1}
 The REK method
\begin{enumerate}[]
\item \mbox{INPUT:} ~$A$, $b$, $\ell$, $x_{0}\in \mathcal{R}(A^T)$ and $z_0\in b+\mathcal{R}(A)$
\item \mbox{OUTPUT:} ~$x_\ell$
\item For $k=1, 2, \ldots, \ell-1$ do
\item ~~~~Select $j\in [n]$ with probability $\frac{\|A_{(j)}\|^2_2}{\|A\|^2_F}$
\item ~~~~Set $z_{k}=z_{k-1}-\frac{A_{(j)}^T z_{k-1}}{ \| A_{(j)} \|_{2}^{2}}A_{(j)}$
\item ~~~~Select $i\in [m]$ with probability $\frac{\|A^{(i)}\|^2_2}{\|A\|^2_F}$
\item ~~~~Set $x_{k}=x_{k-1}-\frac{A^{(i)} x_{k-1}-b^{(i)}+z_{k}^{(i)}}{ \| A^{(i)} \|_{2}^{2}}(A^{(i)})^T$
\item End for
\end{enumerate}
\end{alg}

In \cite{ Dukui2019}, Du presented a tight upper bound for the convergence of the REK method: 
\begin{align}
\mathbb{E}[\|x_{k}-x_\star\|^{2}_2]
\leq \frac{k}{\|A\|^2_F}(1-\frac{\sigma_r^2(A)}{\|A\|^2_F})^k\| z_{0}-(b-Ax_\star )\|^{2}_2+(1-\frac{\sigma_r^2(A)}{\|A\|^2_F})^k\| x_{0}-x_\star \|^{2}_2. \label{1eq}
\end{align}

\section{Convergence analysis}\label{sec3}
A lemma is first given as follows, which will be used to analyze the convergence property of the REK method.
\begin{lem}
\label{theorem1}
Let $A\in R^{m\times n}$, $b\in R^{m}$, $x_{\star}=A^{\dag}b$ be the minimum $\ell_2 $-norm least squares solution, and $v_{\ell}$ be a right singular vector corresponding to the singular value $\sigma_\ell(A)$ of $A$. Let $z_k$ be the $k$th approximation of the RK method applied to $A^Tz=0$ with initial guess $z_0\in b+\mathcal{R}(A)$. Then
\begin{equation}
\label{lem11}
\mathbb{E}[\langle z_{k}-(b-Ax_{\star}), Av_{\ell} \rangle]=(1-\frac{\sigma_\ell^2(A)}{\|A\|^2_F})^k \langle z_{0}-(b-Ax_{\star}), Av_{\ell} \rangle.
\end{equation}
\end{lem}
 \begin{pf}
Let $\mathbb{E}_{k-1}[\cdot]$ be the conditional expectation conditioned on the first $k-1$ iterations of the RK method. Then, from Algorithm \ref{alg1}, we have
\begin{align}
 & \mathbb{E}_{k-1}[\langle z_{k}-(b-Ax_{\star}), Av_{\ell} \rangle] \notag
\\
 &= \mathbb{E}_{k-1}[\langle z_{k-1}-\frac{A_{(j)}^T z_{k-1}}{ \| A_{(j)} \|_{2}^{2}}A_{(j)}-(b-Ax_{\star}), Av_{\ell} \rangle ] \notag
\\
 &= \langle z_{k-1}-(b-Ax_{\star}), Av_{\ell} \rangle -\sum\limits_{j=1}^{n}\frac{ \|A_{ (j )} \|_{2}^{2}}{\|A\|_{F}^{2}}\langle \frac{A_{(j)}^T z_{k-1}}{ \| A_{(j)} \|_{2}^{2}}A_{(j)}, Av_{\ell} \rangle \notag
\\
 &= \langle z_{k-1}-(b-Ax_{\star}), Av_{\ell} \rangle - \frac{ 1}{\|A\|_{F}^{2}}\sum\limits_{j=1}^{n}\langle A_{(j)}^T z_{k-1} A_{(j)}, Av_{\ell} \rangle \notag
\\
 &= \langle z_{k-1}-(b-Ax_{\star}), Av_{\ell} \rangle - \frac{ 1}{\|A\|_{F}^{2}}\sum\limits_{j=1}^{n}\langle A_{(j)}, z_{k-1}  \rangle \langle A_{(j)}, Av_{\ell} \rangle \notag
\\
 &= \langle z_{k-1}-(b-Ax_{\star}), Av_{\ell} \rangle - \frac{ 1}{\|A\|_{F}^{2}} \langle A^Tz_{k-1}, A^TAv_{\ell} \rangle.\notag
\end{align}
Further, by making use of $z_{k-1}=\sum\limits_{i=1}^{m}\langle z_{k-1}, u_i\rangle u_i$, $A^TAv_{\ell}=\sigma_{\ell}^2(A) v_{\ell}$ and $A^Tu_{\ell}=\sigma_{\ell}(A) v_{\ell} $, we get
\begin{align}
 & \mathbb{E}_{k-1}[\langle z_{k}-(b-Ax_{\star}), Av_{\ell} \rangle] \notag
\\
 &= \langle z_{k-1}-(b-Ax_{\star}), Av_{\ell} \rangle - \frac{ 1}{\|A\|_{F}^{2}} \langle A^T\sum\limits_{i=1}^{m}\langle z_{k-1}, u_i\rangle u_i, \sigma_{\ell}^2(A)v_{\ell} \rangle \notag
\\
 &= \langle z_{k-1}-(b-Ax_{\star}), Av_{\ell} \rangle - \frac{ 1}{\|A\|_{F}^{2}} \langle  \sum\limits_{i=1}^{m}\langle z_{k-1}, u_i\rangle \sigma_i(A) v_i, \sigma_{\ell}^2(A)v_{\ell} \rangle, \notag
\end{align}
which together with the orthogonality of the right singular vectors $v_i$ and the fact $Av_i=\sigma_i(A) u_i $ yields
\begin{align}
 & \mathbb{E}_{k-1}[\langle z_{k}-(b-Ax_{\star}), Av_{\ell} \rangle ]\notag
\\
 &= \langle z_{k-1}-(b-Ax_{\star}), Av_{\ell} \rangle - \frac{\sigma_{\ell}^2(A)}{\|A\|_{F}^{2}} \langle  z_{k-1}, u_{\ell} \rangle \sigma_{\ell} (A) \notag
\\
 &= \langle z_{k-1}-(b-Ax_{\star}), Av_{\ell} \rangle - \frac{\sigma_{\ell}^2(A)}{\|A\|_{F}^{2}} \langle  z_{k-1},  \sigma_{\ell} (A)u_{\ell} \rangle \notag
\\
 &= \langle z_{k-1}-(b-Ax_{\star}), Av_{\ell} \rangle - \frac{\sigma_{\ell}^2(A)}{\|A\|_{F}^{2}} \langle  z_{k-1},  Av_{\ell} \rangle. \notag
\end{align}
Thus, by taking the full expectation on both sides and using the facts $A^T(b-Ax_{\star})=0$ and $\langle z_{k}, Av_{\ell} \rangle=\langle A^Tz_{k},  v_{\ell} \rangle $, we have
\begin{align}
 & \mathbb{E}[ \langle z_{k}-(b-Ax_{\star}), Av_{\ell} \rangle ]\notag
\\
 &= \mathbb{E}[\langle z_{k-1}-(b-Ax_{\star}), Av_{\ell} \rangle ]- \frac{\sigma_{\ell}^2(A)}{\|A\|_{F}^{2}} \mathbb{E}[\langle  z_{k-1},  Av_{\ell} \rangle ] \notag
\\
 &= \mathbb{E}[\langle z_{k-1}-(b-Ax_{\star}), Av_{\ell} \rangle ]- \frac{\sigma_{\ell}^2(A)}{\|A\|_{F}^{2}} \mathbb{E}[\langle A^T (z_{k-1}-(b-Ax_{\star})),   v_{\ell} \rangle ] \notag
\\
 &= \mathbb{E}[\langle z_{k-1}-(b-Ax_{\star}), Av_{\ell} \rangle] - \frac{\sigma_{\ell}^2(A)}{\|A\|_{F}^{2}} \mathbb{E}[\langle  (z_{k-1}-(b-Ax_{\star})),   Av_{\ell} \rangle]  \notag
\\
 &= (1- \frac{\sigma_{\ell}^2(A)}{\|A\|_{F}^{2}} )\mathbb{E}[\langle  (z_{k-1}-(b-Ax_{\star})),   Av_{\ell} \rangle].  \notag
\end{align}
By induction on the iteration index $k$, we can obtain the estimate (\ref{lem11}).
 \end{pf}

In the following, we give two main observations of the REK method.

\begin{thm}
\label{theorem2}
Let $A\in R^{m\times n}$, $b\in R^{m}$, $x_{\star}=A^{\dag}b$ be the minimum $\ell_2 $-norm least squares solution, and $v_{\ell}$ be a right singular vector corresponding to the singular value $\sigma_\ell(A)$ of $A$.
Let $x_k$ be the $k$th 
approximation of the REK method generated by Algorithm \ref{alg1} with initial guess $x_{0}\in \mathcal{R}(A^T)$ and $z_0\in b+\mathcal{R}(A)$. Then
\begin{align}
\mathbb{E}[\langle x_{k}- x_{\star}, v_{\ell} \rangle]=\frac{k}{\|A\|^2_F}(1-\frac{\sigma_\ell^2(A)}{\|A\|^2_F})^k \langle -A^Tz_{0}, v_{\ell} \rangle+ (1-\frac{\sigma_\ell^2(A)}{\|A\|^2_F})^k\langle x_{0}- x_{\star}, v_{\ell} \rangle. \label{th}
\end{align}
\end{thm}

\begin{pf}
Since
\begin{align}
 \mathbb{E}[ \langle x_{k}- x_{\star}, v_{\ell} \rangle ] = \mathbb{E} [\langle x_{k}- \hat{x}_{k}, v_{\ell} \rangle]+\mathbb{E} [\langle\hat{x}_{k}- x_{\star}, v_{\ell} \rangle], \label{eq5}
\end{align}
where $\hat{x}_{k}$ is the one-step 
update of the RK method for solving $Ax=Ax_{\star}$ from $x_{k-1}$, i.e., $\hat{x}_{k}=x_{k-1}-\frac{A^{(i)}x_{k-1}-A^{(i)}x_{\star}}{\|A^{(i)}\|^2_2}(A^{(i)})^T$, we next 
consider $\mathbb{E}[ \langle x_{k}- \hat{x}_{k}, v_{\ell} \rangle]$ and $\mathbb{E} [\langle\hat{x}_{k}- x_{\star}, v_{\ell} \rangle]$ separately.

We first consider $ \mathbb{E} [\langle x_{k}- \hat{x}_{k}, v_{\ell} \rangle]$. Let $\mathbb{E}_{k-1}[\cdot]$ be the conditional expectation conditioned on the first $k-1$ iterations of the REK method. That is,
\begin{align}
 \mathbb{E}_{k-1}[\cdot]= \mathbb{E}[\cdot|j_1, i_1, j_2, i_2, \ldots, j_{k-1}, i_{k-1}],  \notag
\end{align}
where $j_{t^*}$ is the ${t^*}$th column chosen and $i_{t^*}$ is the ${t^*}$th row chosen. We denote the conditional expectation conditioned on
the first $k-1$ iterations and the $k$th column chosen as
\begin{align}
 \mathbb{E}_{k-1}^{i}[\cdot]= \mathbb{E}[\cdot|j_1, i_1, j_2, i_2, \ldots, j_{k-1}, i_{k-1}, j_k]. \notag
\end{align}
Similarly, we denote the conditional expectation conditioned on the first $k-1$ iterations and the $k$th row chosen as
\begin{align}
 \mathbb{E}_{k-1}^{j}[\cdot]= \mathbb{E}[\cdot|j_1, i_1, j_2, i_2, \ldots, j_{k-1}, i_{k-1}, i_k]. \notag
\end{align}
Then, by the law of total expectation, we have
\begin{align}
 \mathbb{E}_{k-1}[\cdot]=  \mathbb{E}_{k-1}^{j}[ \mathbb{E}_{k-1}^{i}[\cdot] ]. \notag
\end{align}
Thus, according to the update formulas of $\hat{x}_{k}$ given above and $x_k$ given in Algorithm \ref{alg1}, 
we obtain
\begin{align}
 &\mathbb{E}_{k-1} [\langle x_{k}- \hat{x}_{k}, v_{\ell} \rangle ]\notag
\\
 &= \mathbb{E}_{k-1}[\langle x_{k-1}-\frac{A^{(i)} x_{k-1}-b^{(i)}+z_{k}^{(i)}}{ \| A^{(i)} \|_{2}^{2}}(A^{(i)})^T- (x_{k-1}-\frac{A^{(i)}x_{k-1}-A^{(i)}x_{\star}}{\|A^{(i)}\|^2_2}(A^{(i)})^T), v_{\ell} \rangle] \notag
\\
 &= \mathbb{E}_{k-1}[\langle \frac{b^{(i)}-A^{(i)} x_{\star}-z_{k}^{(i)}}{ \| A^{(i)} \|_{2}^{2}}(A^{(i)})^T, v_{\ell} \rangle ]\notag
\\
 &=\mathbb{E}_{k-1}^{j}[ \mathbb{E}_{k-1}^{i} [\langle \frac{b^{(i)}-A^{(i)} x_{\star}-z_{k}^{(i)}}{ \| A^{(i)} \|_{2}^{2}}(A^{(i)})^T, v_{\ell} \rangle ] ] \notag
\\
 &=\mathbb{E}_{k-1}^{j}[ \frac{1}{\|A\|^2_F} \sum\limits_{i=1}^{m} \langle (b^{(i)}-A^{(i)} x_{\star}-z_{k}^{(i)})(A^{(i)})^T, v_{\ell} \rangle  ] \notag
\\
 &=\mathbb{E}_{k-1}^{j}[ \frac{1}{\|A\|^2_F} \sum\limits_{i=1}^{m} (b^{(i)}-A^{(i)} x_{\star}-z_{k}^{(i)})\langle (A^{(i)})^T, v_{\ell} \rangle  ] \notag
\\
 &=\mathbb{E}_{k-1}^{j}[ \frac{1}{\|A\|^2_F}  \langle b -A  x_{\star}-z_{k} , Av_{\ell} \rangle  ] \notag
\\
 &=\frac{1}{\|A\|^2_F}  \mathbb{E}_{k-1} [  \langle b -A  x_{\star}-z_{k} , Av_{\ell} \rangle  ].\notag
\end{align}
As a result,
\begin{align}
 &\mathbb{E} [ \langle x_{k}- \hat{x}_{k}, v_{\ell} \rangle] =\frac{1}{\|A\|^2_F}  \mathbb{E} [  \langle b -A  x_{\star}-z_{k}, Av_{\ell} \rangle  ], \notag
\end{align}
which together with Lemma \ref{theorem1} yeilds
\begin{align}
 &\mathbb{E}  [\langle x_{k}- \hat{x}_{k}, v_{\ell} \rangle ]=\frac{1}{\|A\|^2_F}  (1-\frac{\sigma_\ell^2(A)}{\|A\|^2_F})^k \langle b-Ax_{\star}- z_{0}, Av_{\ell}  \rangle. \label{eq6}
\end{align}

We now consider $\mathbb{E} [\langle\hat{x}_{k}- x_{\star}, v_{\ell} \rangle]$. It follows from $\hat{x}_{k}=x_{k-1}-\frac{A^{(i)}x_{k-1}-A^{(i)}x_{\star}}{\|A^{(i)}\|^2_2}(A^{(i)})^T$ that
\begin{align}
 &\mathbb{E}_{k-1}  [\langle\hat{x}_{k}- x_{\star}, v_{\ell} \rangle ]\notag
\\
 &= \mathbb{E}_{k-1} [ \langle x_{k-1}-\frac{A^{(i)}x_{k-1}-A^{(i)}x_{\star}}{\|A^{(i)}\|^2_2}(A^{(i)})^T- x_{\star}, v_{\ell} \rangle ]  \notag
\\
 &=  \langle x_{k-1}-x_{\star}, v_{\ell} \rangle -   \frac{1}{\|A\|^2_F} \sum\limits_{i=1}^{m} \langle   A^{(i)}(x_{k-1}- x_{\star})(A^{(i)})^T, v_{\ell} \rangle \notag
\\
 &=  \langle x_{k-1}-x_{\star}, v_{\ell} \rangle -   \frac{1}{\|A\|^2_F} \sum\limits_{i=1}^{m}\langle   (A^{(i)})^T, x_{k-1}- x_{\star} \rangle \langle  (A^{(i)})^T, v_{\ell} \rangle \notag
\\
 &=  \langle x_{k-1}-x_{\star}, v_{\ell} \rangle -   \frac{1}{\|A\|^2_F}  \langle   A( x_{k-1}- x_{\star}), A v_{\ell}\rangle  \notag
\\
 &=  \langle x_{k-1}-x_{\star}, v_{\ell} \rangle -   \frac{1}{\|A\|^2_F}  \langle   A(\sum\limits_{i=1}^{n} \langle  ( x_{k-1}- x_{\star}), v_i\rangle v_i), A v_{\ell}\rangle  \notag
\\
 &=  \langle x_{k-1}-x_{\star}, v_{\ell} \rangle -   \frac{1}{\|A\|^2_F}  \langle   \sum\limits_{i=1}^{n} \langle  ( x_{k-1}- x_{\star}), v_i\rangle \sigma_i(A)u_i , \sigma_{\ell}(A)u_{\ell}\rangle  \notag
\\
 &=  \langle x_{k-1}-x_{\star}, v_{\ell} \rangle -   \frac{\sigma_{\ell}^2(A)}{\|A\|^2_F}   \langle x_{k-1}-x_{\star}, v_{\ell} \rangle \notag
\\
 &=  (1 -   \frac{\sigma_{\ell}^2(A)}{\|A\|^2_F})   \langle x_{k-1}-x_{\star}, v_{\ell} \rangle. \notag
\end{align}
Thus
\begin{align}
 \mathbb{E}   [\langle\hat{x}_{k}- x_{\star}, v_{\ell} \rangle ] =  (1 -   \frac{\sigma_{\ell}^2(A)}{\|A\|^2_F})  \mathbb{E}   [\langle x_{k-1}-x_{\star}, v_{\ell} \rangle]. \label{eq7}
\end{align}

Combining (\ref{eq5}), (\ref{eq6}), and (\ref{eq7}) yields
\begin{align}
\mathbb{E} [\langle x_{k}- x_{\star}, v_{\ell} \rangle ]&= \mathbb{E}[ \langle x_{k}- \hat{x}_{k}, v_{\ell} \rangle]+\mathbb{E} [\langle\hat{x}_{k}- x_{\star}, v_{\ell} \rangle ]\notag
\\
&=\frac{1}{\|A\|^2_F}  (1-\frac{\sigma_\ell^2(A)}{\|A\|^2_F})^k \langle b-Ax_{\star}- z_{0}, Av_{\ell}  \rangle + (1 -   \frac{\sigma_{\ell}^2(A)}{\|A\|^2_F})  \mathbb{E}  [ \langle x_{k-1}-x_{\star}, v_{\ell} \rangle ]\notag
\\
&=\frac{k}{\|A\|^2_F}  (1-\frac{\sigma_\ell^2(A)}{\|A\|^2_F})^k \langle b-Ax_{\star}- z_{0}, Av_{\ell}  \rangle + (1 -   \frac{\sigma_{\ell}^2(A)}{\|A\|^2_F})^k  \langle x_{0}-x_{\star}, v_{\ell} \rangle \notag
\\
&=\frac{k}{\|A\|^2_F}  (1-\frac{\sigma_\ell^2(A)}{\|A\|^2_F})^k \langle  - A^Tz_{0},  v_{\ell}  \rangle + (1 -   \frac{\sigma_{\ell}^2(A)}{\|A\|^2_F})^k  \langle x_{0}-x_{\star}, v_{\ell} \rangle, \notag
\end{align}
which implies the desired result (\ref{th}).
\end{pf}

\begin{rmk}
\label{rmk1}
Theorem \ref{theorem2} shows that the decay rates of $\|x_k-x_{\star}\|_2$ are different in different right singular vectors spaces. Specifically, the decay rates of the REK method are dependent on the singular values: the larger the singular value of 
$A$ is, the faster the error decays in the corresponding right singular vector space. This implies that the smallest singular value will lead to the slowest rate of convergence, which is the 
one in (\ref{1eq}). So, 
the convergence bound presented by Du \cite{ Dukui2019} is optimal. The above findings also explain the phenomenon found
in the extensive numerical experiments appearing in the literature that the REK method seems to converge faster in the beginning.
\end{rmk}

\begin{rmk}
\label{rmk11}
If the system (\ref{1}) is consistent, i.e., $b\in \mathcal{R}(A)$, which implies that $z_0\in\mathcal{R}(A)$, 
Theorem \ref{theorem2} automatically reduces to the following result: 
\begin{align}
\mathbb{E} [\langle x_{k}- x_{\star}, v_{\ell} \rangle] =  (1 -   \frac{\sigma_{\ell}^2(A)}{\|A\|^2_F})^k  \langle x_{0}-x_{\star}, v_{\ell} \rangle, \notag
\end{align}
which recovers the decay rates of the RK method in different right singular vectors spaces \cite{steinerberger2021randomized}.
\end{rmk}
\begin{thm}
\label{theorem3}
Let $A\in R^{m\times n}$, $b\in R^{m}$ and $x_{\star}=A^{\dag}b$ be the minimum $\ell_2 $-norm least squares solution. Let $x_k$ be the $k$th 
approximation of the REK method generated by Algorithm \ref{alg1} with initial guess $x_{0}\in \mathcal{R}(A^T)$ and $z_0\in b+\mathcal{R}(A)$. Then
\begin{align}
\mathbb{E}[ \|x_{k}- x_{\star}\|^2_2]\leq\frac{1}{\|A\|^2_F}(1-\frac{\sigma_r^2(A)}{\|A\|^2_F})^k \| z_{0}-(b-Ax_{\star})\|^2_2+ \mathbb{E}[(1-\frac{1}{\|A\|^2_F}\|A\frac{x_{k-1}-x_{\star}}{\|x_{k-1}-x_{\star}\|_2}\|_2^2)\|x_{k-1}-x_{\star}\|_2^2]. \label{th2}
\end{align}
\end{thm}

\begin{pf}
Following an analogous argument to Theorem 2 
of \cite{Dukui2019}, we get
\begin{align}
\mathbb{E}   [\|x_{k}- x_{\star}\|_2^2 ]&= \mathbb{E}[ \| x_{k}- \hat{x}_{k}\|_2^2]+\mathbb{E} [\|\hat{x}_{k}- x_{\star}\|_2^2], \notag
\end{align}
\begin{align}
\mathbb{E} [\| x_{k}- \hat{x}_{k}\|_2^2]\leq\frac{1}{\|A\|^2_F}(1-\frac{\sigma_r^2(A)}{\|A\|^2_F})^k \| z_{0}-(b-Ax_{\star})\|^2_2, \notag
\end{align}
and
\begin{align}
\mathbb{E} [\|\hat{x}_{k}- x_{\star}\|_2^2]
&= \mathbb{E} [(x_{k-1}-x_{\star})^T(I-\frac{A^TA}{\|A\|_F^2})(x_{k-1}-x_{\star})]\notag
\\
&=\mathbb{E}[(\|x_{k-1}-x_{\star}\|_2^2- \frac{1}{\|A\|^2_F}\|A(x_{k-1}-x_{\star})\|_2^2)]\notag
\\
&=\mathbb{E}[(1-\frac{1}{\|A\|^2_F}\|A\frac{x_{k-1}-x_{\star}}{\|x_{k-1}-x_{\star}\|_2}\|_2^2)\|x_{k-1}-x_{\star}\|_2^2]. \notag
\end{align}
Combining the above three equations, we have
\begin{align}
\mathbb{E}  [ \|x_{k}- x_{\star}\|_2^2]
&= \mathbb{E} [\| x_{k}- \hat{x}_{k}\|_2^2]+\mathbb{E}[ \|\hat{x}_{k}- x_{\star}\|_2^2 ] \notag
\\
&\leq\frac{1}{\|A\|^2_F}(1-\frac{\sigma_r^2(A)}{\|A\|^2_F})^k \| z_{0}-(b-Ax_{\star})\|^2_2+ \mathbb{E}[(1-\frac{1}{\|A\|^2_F}\|A\frac{x_{k-1}-x_{\star}}{\|x_{k-1}-x_{\star}\|_2}\|_2^2)\|x_{k-1}-x_{\star}\|_2^2], \notag
\end{align}
which implies the desired result (\ref{th2}).
\end{pf}

\begin{rmk}
\label{rmk2}
Since $\|A\frac{x_{k-1}-x_{\star}}{\|x_{k-1}-x_{\star}\|_2}\|_2^2\geq\sigma_r^2(A)$, Theorem \ref{theorem3} 
implies that the REK method actually converges faster if $x_{k-1}-x_{\star}$ is not close to right singular vectors corresponding to the small singular values of $A$ .

\end{rmk}
\section{Numerical experiment}\label{sec4}
\begin{figure}[ht]
 \begin{center}
\includegraphics [height=5.5cm,width=8.5cm ]{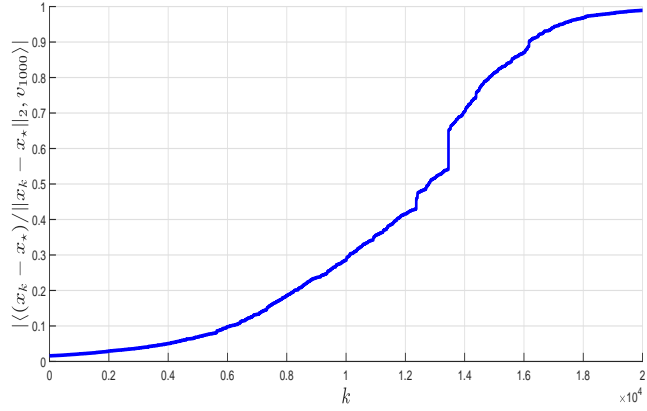}
 \end{center}
\caption{A sample evolution of $  |\langle (x_k-x_{\star})/\|x_k-x_{\star}\|_2, v_{1000} \rangle|$.}\label{fig1}
\end{figure}

\begin{figure}[ht]
 \begin{center}
\includegraphics [height=5.5cm,width=8.5cm  ]{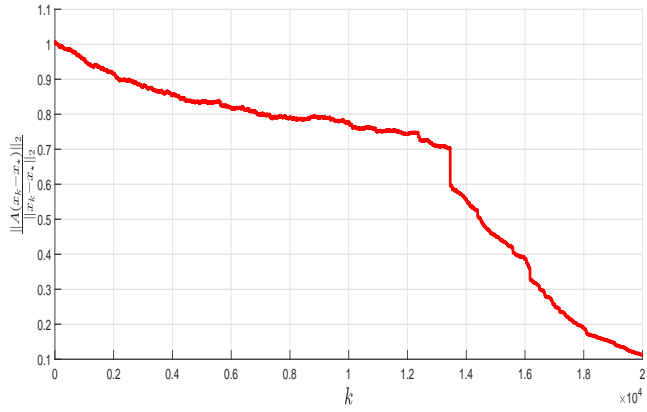}
 \end{center}
\caption{A sample evolution of $\frac{\|A(x_k-x_{\star})\|_2 }{\|x_k-x_{\star}\|_2}$. }\label{fig2}
\end{figure}

Now we present a simple example to illustrate that the REK method finally converges along right singular vector corresponding to the small singular value of $A$. 
To this end, let $A_0\in R^{1000\times 1000}$ be a Gaussian matrix with i.i.d. $N(0, 1)$ entries and $D\in R^{1000\times 1000}$ be a diagonal matrix whose diagonal elements are all 100. Further, we set $A_1=A_0+D$ and replace 
its last row $A_1^{(1000)}$ by a tiny perturbation of $A_1^{(999)}$, i.e., 
adding 0.01 to each entry of $A_1^{(999)}$. Then, we normalize all rows of $A_1$, i.e., set $\|A_1^{(i)}\|_2=1$, $i=1, 2, \ldots, 1000$. After that, we set $A=\begin{bmatrix}
A_1\\
A_2
\end{bmatrix}
\in R^{1100\times 1000}$, where $A_2\in R^{100\times 1000}$ is a zero matrix. So, the first 999 singular values of the matrix $A$ are between $\sim 0.5$ and $\sim 1.5$, and the smallest nonzero singular value is $\sim 10^{-4}$. In addition, we generate the solution vector $x_\star\in R^{1000}$ using the MATLAB function \texttt{randn}, and set the right-hand side $b=Ax_{\star}+z$, where $z$ is a nonzero vector belonging to the null space of $A^{T}$, which is generated by the MATLAB function \texttt{null}. That is, the system (\ref{1}) is inconsistent. With $x_0=0$ and $z_0=b$, we plot $|\langle (x_k-x_{\star})/\|x_k-x_{\star}\|_2, v_{1000} \rangle|$ and $\frac{\|A(x_k-x_{\star})\|_2 }{\|x_k-x_{\star}\|_2}$ in Figure \ref{fig1} and Figure \ref{fig2}, respectively.

From Figure \ref{fig1}, we find that $|\langle (x_k-x_{\star})/\|x_k-x_{\star}\|_2, v_{1000} \rangle|$ initially is very small and almost is 0, which indicates that $x_k-x_{\star} $ is not 
close to the right singular vector $v_{1000}$. Considering the analysis of Remark \ref{rmk2}, the phenomenon implies the `preconvergence' behavior of the REK method, that is, the REK method seems to converge quickly at the beginning. In addition, as $k\rightarrow\infty$, $|\langle (x_k-x_{\star})/\|x_k-x_{\star}\|_2, v_{1000} \rangle|\rightarrow 1$. This phenomenon implies that $x_{k}-x_{\star}$ tends to the right singular vector corresponding to the smallest singular value of $A$.

From Figure \ref{fig2}, we observe that as $k$ increases, $\frac{\|A(x_k-x_{\star})\|_2 }{\|x_k-x_{\star}\|_2}$ approaches the small singular value.
This phenomenon implies 
the same result 
given above, i.e., 
as $k\rightarrow\infty$, $x_{k}-x_{\star}$ tends to the right singular vector corresponding to the smallest singular value of $A$.
Furthermore, this phenomenon also allows for an interesting application, i.e., finding nonzero vectors $x$ such that $\frac{\|Ax\|_2}{\|x\|_2}$ is small.

\clearpage

\bibliography{mybibfile}

\end{document}